\documentclass{article}
\usepackage{amsmath}
\usepackage{amssymb}

\title{Robinson-Schensted correspondence and left cells}
\author{Susumu Ariki}
\date{}

\newtheorem{thm}{Theorem}[section]
\newtheorem{example}[thm]{Example}
\newtheorem{defn}[thm]{Definition}
\newtheorem{prop}[thm]{Proposition}
\newtheorem{lemma}[thm]{Lemma}

\newcommand{\qed}{\hbox{Q.E.D}}

\begin{document}
\maketitle

\section{Introduction}
This is based on \cite{A}. 
In \cite{A}, I explained several 
theorems which focused on a famous theorem 
of \cite{KL} that 
two elements of the symmetric group belong to 
a same left cell if and only if they share 
a common Q-symbol. The first half of \cite{A} was 
about the direct proof of this theorem (Theorem A), and 
the second half was about the relation between 
primitive ideals and left cells, and I explained 
another proof of this theorem. 

The reason why I gave a direct 
proof which was different from the proof in 
\cite{KL} was that the proof in \cite{KL} was 
hard to read: 
It relied on \cite[6]{V1}, which in turn 
relied on \cite{Jo1}, the full paper of 
which is not yet available even today. Note also that 
the theorem itself is not stated in \cite{KL}. But 
the beginning part of the proof of \cite[Theorem 1.4]{KL} 
gives some explanation on the relation between 
left cells in the sense of 
Kazhdan and Lusztig and Vogan's generalized 
$\tau$-invariants in the theory of 
primitive ideals. In this picture, 
Theorem A is derived from Joseph's theorem. 

Lack of a clear proof in the literature lead 
Garsia and McLarnan to the publication of \cite{GM}.
\footnote{It is worth mentioning that 
Garsia and McLarnan wrote in \cite{GM} that they were benefitted by 
A.Bj\"orner's lecture notes and R.King's 
lecture notes. Both of these notes 
are still not available, and it seems that 
preliminary version of them were circulated 
in a very restricted group of people around the time. }
The proof given in \cite{GM} is close to \cite{A}, 
but the line of the proof 
in \cite{GM} is interrupted with 
combinatorics of tableaux, which is not 
necessary. In fact, after we read to the fourth section of 
\cite{KL}, which is the section for some preliminaries 
to the proof of \cite[Theorem 1.4]{KL}, we can give a 
short and elementary proof of Theorem A in a direct way, 
as I will show below. 

I rush to say that my proof was not so original: 
It copied argument in \cite[Satz 5.25]{Ja1} for 
Joseph's theorem. This was the reason why I did not publish 
it in English. But after a decade has passed, we still have 
no suitable literature which includes the direct proof. 
Further, we have new development in the last decade. 
For example, we have better understanding of this theorem in 
{\it jeu de taquin} context (see \cite{H},\cite{BSS}); 
the study of sovable lattice models in 
Kyoto school lead to the theory of crystal bases and 
canonical bases, by which we can understand Theorem A 
in the crystal base theory context. 

I have therefore decided to add this 
short note to this volume in order to give this proof and 
new development. 
I also prepare enough papers 
in the references for reader's convenience. 
I note that there is a sketch of a proof of Theorem A in \cite[p.172]{BV}. 
It involves the notion of wave front sets. 
I do not recommend \cite{BV} for knowing the proof of Theorem A. 
One reason is that we do not need wave front sets for the proof 
of Theorem A itself. 

The non-direct proof I explained in the second half of \cite{A} 
is the proof which is indicated in \cite{KL}, and the proof was taken from 
Jantzen's lecture notes \cite{Ja1}. Thus I do not reproduce it, 
and I only give statements 
of several theorems (Theorem B,C,D) which concludes Theorem A. 

I give some bibliographical comments on the second proof. 
It is obtained by combining two theorems; 
one known as Joseph's theorem, 
which states that two primitive ideals with a same 
integral regular central character coincide if and only if 
their Q-symbol coincide (Theorem B), and another theorem due to 
Joseph and Vogan, which relates the inclusion relation 
of primitive ideals to non-vanishing condition 
of certain multiplicities (Theorem C), and thus to 
order relation of left cells of 
the symmetric group (Theorem D). 

In a survey \cite{Bo}, it is stated 
that Theorem A was proved in \cite{Jo1}, and simple 
proof could be found in \cite{V1} and \cite{Ja1}. But 
as I stated, \cite{Jo1} is not published, and \cite{V1} is 
based on \cite{Jo1}. Thus to read \cite{V1}, 
one has to reproduce arguments by oneself. 
Nevertheless, it is Joseph's theorem, and his idea came from 
the explicit form of Goldie rank representations in type A case 
\cite[Proposition 8.4]{Jo2}, which makes 
the number of primitive ideals with a common regular central character 
equal to the number of involutions. That is, 
Duflo's map is bijective, which proves Theorem B. 
The proof of Theorem B 
is given in \cite[Corollary 5.3]{Jo3}. We do not follow his 
line and I refer to \cite[Satz 5.25]{Ja1}. 

Theorem C is proved in \cite[Theorem 3.2]{V2}. One 
implication is due to \cite[Theorem 5.3]{Jo5}, which is 
reformulated in \cite[Proposition 3.1]{V2}. 
It is not difficult to derive Theorem D: That 
Theorem C implies Theorem D (see also \cite[Conjecture C]{Jo5}) is 
stated in the introduction of \cite{V2}. The proof here is 
based on \cite[Corollar 7.13]{Ja1} and \cite[Lemma 14.9]{Ja1}. 
Since the Kazhdan-Lusztig conjecture proved by 
Brylinski-Kashiwara and Beilinson-Berstein is very well known, 
I take it for granted when I explain Theorem D. 
But it is of cource a very deep result. 

Joseph's Goldie rank representation was related to Springer 
representation and the theory evolved into a beautiful 
geometric representation theory. Since this part is not at all 
combinatorial, I do not mention it. 
I only refer to \cite{BB3} for this development. 
There is also new direction for the genaralization of the 
Robinson-Schensted correspondence related to the primitive 
ideal theory. I refer \cite{Ga} and \cite{Tr}. 
\footnote{There is one more dirction: 
generalization of the Steinberg's theorem 
\cite{St} is given by M. van Leeuwen. This is the direction 
to the geometry of flag varieties.}

\section{Preliminaries}
\subsection{P-symbols and Q-symbols}

Let $S_n$ be the symmetric group of degree $n$. 
Namely, its underlying set is the set of 
bijective maps from $\{1,\dots,n\}$ to itself, 
and the group structure is given by 
composition of maps. Let $w$ be an element in $S_n$, 
and denote 
the image of $i\in\{1,\dots,n\}$ under the map $w$ 
by $w_i$. Throughout the paper, we identify $w$ 
with the sequence $w_1\cdots w_n$, 
which is a permutation of $1,\dots,n$. 

Let ${\mathbb N}=\{1,2,\dots\}$ be the set of 
natural numbers. 
A {\bf Young diagram} $\lambda$ is a finite subset of 
${\mathbb N}\times{\mathbb N}$ which satisfies 
the condition that if $(x,y)\in\lambda$, then 
$\{(x\!-\!1,y),(x,y\!-\!1)\}\cap
{\mathbb N}\times{\mathbb N}\in\lambda$. 
$(x,y)\in\lambda$ is called a {\bf cell} of $\lambda$. 
$x$ is called the row number of the cell, and 
$y$ is called the column number of the cell. 

A {\bf tableau} $T$ of shape $\lambda$ is a map 
from $\lambda$ to $\mathbb N$. 
The image of a cell $(x,y)$ of $\lambda$ under the map is called 
the entry of the cell, and is denoted by $T(x,y)$. We only consider 
the tableaux satisfying 
\[
T(x,y)\le T(x',y')\quad(x\le x',y\le y'). 
\]
If it also satisfies $T(x,y)<T(x',y)\;(x\!<\!x')$ 
(resp. $T(x,y)<T(x,y')\;(y\!<\!y'))$, $T$ is 
called a {\bf column strict} 
(resp. {\bf row strict}) {\bf semi-standard tableau}. 
If the entries of $T$ are precisely $\{1,\dots,n\}$, 
we call $T$ a {\bf standard tableau}. 

Let $T$ be a column strict semi-standard tableau, 
$k$ be a natural number. 
We denote the set of cells in the $i$ th row 
by $R_i(T)$, and denote the maximal column number of the 
cells in 
$R_i(T)$ by $c_i$. Assume that we are given a natural 
number $k_i$. If $k_i$ is equal or greater than all entries 
of $R_i(T)$, we add the cell $(i,c_i\!+\!1)$ to 
$R_i(T)$ and make its entry be $k_i$. 
If it is not the case, we consider the cells 
of $R_i(T)$ whose entries are greater than $k_i$, 
and pick up the cell of minimal column number among them. 
We then change its entry to $k_i$, and we make 
$k_{i+1}$ be the original entry of the cell. 
This latter procedure is called {\bf bumping procedure}. 
We set $k_1=k$, and continue the bumping procedure 
until no bumping occurs. This is called a 
{\bf row insertion algorithm}, and it results in a 
new column strict semi-standard tableau, which we denote by 
$T\leftarrow k$. In the similar way, we can define a 
column insertion algorithm $k\rightarrow T$. 

\begin{defn}
Let $w=w_1\cdots w_n$ be a permutation. 
Two standard tableaux 
$P(w)$ and $Q(w)$ 
defined by 
\[
\begin{array}{rl}
P(w)\!&=\emptyset\!\leftarrow\! w_1\!\leftarrow\cdots
\leftarrow\! w_n,\\
Q(w)\!&=P(w^{-1})
\end{array}
\]
are called the {\bf P-symbol} of $w$ and the {\bf Q-symbol} of 
$w$ respectively. The correspondence between 
$w$ and the pair $(P(w),Q(w))$ is called 
the {\bf Robinson-Schensted correspondence}. 
We often write $P(w)=\emptyset\!\leftarrow w_1\cdots w_n$ for short. 
\end{defn}

It is known that $P(w)=w_1\!\rightarrow\!\cdots\!\rightarrow\! w_k
\!\rightarrow\!\emptyset\!\leftarrow\! w_{k+1}\!\leftarrow\!\cdots
\!\leftarrow\! w_n$ for any $k$. 

\begin{example}
If $w=31524$, then we have 
\[
\begin{array}{cccccccc}
P(w)=&1&2&4&\qquad\qquad Q(w)=&1&3&5\\          
     &3&5& &                  &2&4& \end{array}.
\]
\end{example}

\bigskip
More familiar definition of the Q-symbol is by 
the "recording" tableau, which records the cell added by each 
insertion procedure. It is a well known theorem 
that it coincides with $P(w^{-1})$. 

\bigskip
\noindent
{\bf Remark}
We have a two dimensional pictorial algorithm 
to compute P-symbols and Q-symbols due to S.V.Fomin 
\cite[4.2.4]{Fo3}. In his picture, we know at a glance 
that $Q(w)$ equals $P(w^{-1})$. 

\bigskip
\noindent
{\bf Remark}
If two elements in $S_n$ have a common P-symbol, 
we say that these belong to a same left Knuth class. 
Similarly, if these have a common Q-symbol, 
we say that these belong to a same right Knuth class. 

Although we do not go into the combinatorial 
structures of the Robinson- 
Schensted correspondence, 
it is worth referring to the relation of the 
Robinson-Schensted 
correspondence to the {\bf\it jeu de taquin} sliding 
algorithm. Namely, the 
insertion algorithm consists of {\it jeu de taquin} moves, 
and {\it jeu de taquin} equivalence classes are 
the same as left Knuth classes. 
On the other hand, 
right Knuth classes are the same as 
Haiman's dual equivalence classes. 

To be more precise, let $\lambda_n$ be the staircase 
Young diagram of size $n(n\!-\!1)/2$. Then 
$\lambda_{n+1}/\lambda_n$ consists of $n$ one-cell 
components. The tableaux of shape 
$\lambda_{n+1}/\lambda_n$ are called 
{\bf permutation tableaux}. By reading entries from 
left to right, we identify permutation tableaux 
with permutations of $1,\dots,n$. This is in fact 
true for more general tableaux. We read entries 
of such a tableau row by row from left to right 
starting with the cell on the south-west end and ending up with 
the cell on the north-east end. Then the corresponding 
permutation tableau is {\it jeu de taquin} equivalent 
to the original tableau. 

We take a tableau $T$ of shape $\lambda_n$, and 
compute {\bf switching} of a permutation tableau 
$w$ and $T$ \cite{BSS}. 
Since all tableaux of a same non skew shape 
are dual equivalent \cite[Proposition 4.2]{BSS}, 
the non skew tableau produced by the switching 
is independent of the choice of $T$ 
\cite[Theorem 4.3]{BSS}, and it is {\it jeu de taquin} 
equivalent to the permutation tableau $w$ 
\cite[Theorem 3.1]{BSS}. This is the P-symbol of $w$. 
Further, two permutation tableaux are in a same dual 
equivalence class if and only if they have a 
common Q-symbol \cite[Theorem 2.12]{H}. 
These give the Robinson-Schensted corespondence in the 
{\it jeu de taquin} context. In fact, this view point 
also appears in the crystal base theory \cite{BKK}. 

\subsection{KL polynomials}

Let $q$ be a variable, and let ${\mathcal H}_n$ be 
the Hecke algebra of the symmetric group $S_n$. 
Namely, ${\mathcal H}_n$ is 
the algebra over ${\mathbb Q}(q)$ defined by 
generators $T_1,\dots,T_{n-1}$ and relations 
\[
(T_i-q)(T_i+1)=0, \quad 
T_iT_{i+1}T_i=T_{i+1}T_iT_{i+1},\quad 
T_iT_j=T_jT_i \;(j\ge i+2). 
\]
Let $s_i=(i,i\!+\!1)$ be the transposition of 
$i$ and $i\!+\!1$. We set $T_{s_i}=T_i$. 
For general $w\in S_n$, we find a 
reduced expression $w=s_{i_1}\cdots s_{i_r}$ 
and set $T_w=T_{i_1}\cdots T_{i_r}$. 
The length $r$ of reduced expressions does not 
depend on the reduced expressions, 
which is denoted by $l(w)$. It is also 
well known that $T_w$ does not depend on the 
choice of the reduced expression, and 
$\{T_w|w\in S_n\}$ is a basis of 
${\mathcal H}_n$. If $y$ is obtained by the product 
of a subword of a reduced expression of $w$, 
we write $y\le w$. This order is called 
{\bf Bruhat order}. 

\begin{defn} 
The following two conditions uniquely define 
the polynomials $P_{y,w}(q)\in {\mathbb Z}[q]
\;(y\le w)$, which 
are called {\bf Kazhdan-Lusztig polynomials}: 
\[
\begin{array}{lll}
C_w\!&\!=\!&\sum_{y\le w} (-1)^{l(w)-l(y)}
q^{\frac{l(w)}{2}-l(y)}P_{y,w}(q^{-1})T_y \\
&\!=\!& \sum_{y\le w} (-1)^{l(w)-l(y)}
q^{-\frac{l(w)}{2}+l(y)}P_{y,w}(q)T_{y^{-1}}^{-1}, 
\end{array}
\]
and $P_{w,w}(q)=1, \quad
\mbox{\rm deg} P_{y,w}(q)\le\frac{l(w)-l(y)-1}{2} \quad
(y<w).$
\end{defn}

We call the first property the bar invariance 
property, and the second property the degree 
property. For the definition, we can use the following 
element instead of $C_w$. 
\[
C_w'=q^{-\frac{l(w)}{2}}\sum_{y\le w}P_{y,w}(q)T_y
\]

If $\frac{l(w)-l(y)-1}{2}$ is an integer, we denote 
the coefficient of $q^{\frac{l(w)-l(y)-1}{2}}$ 
in $P_{y,w}(q)$ by $\mu(y,w)$. If 
$\mu(y,w)\ne 0$ for $y<w$ or 
$\mu(w,y)\ne 0$ for $y>w$ occurs, 
we write $\mu(y|w)\ne 0$. Note that if we write 
$\mu(y|w)\ne 0$, it particularly implies 
that $y>w$ or $y<w$ holds. 

\bigskip
The welldefinedness of $P_{y,w}(q)$ is non trivial, and 
in fact is one of the main theorems \cite[Theorem 1.1]{KL}. 
The uniqueness of $P_{y,w}(q)$ is easy to prove, but 
for the existence of these polynomials, 
we need to construct $C_w\;(w\in S_n)$. In \cite[2.2]{KL}, 
these $C_w$ are inductively constructed by setting 
$C_e=1$ and 
\[
C_w=C_{s_i}C_{s_iw}-\underset{\mu(z,w)\ne0}
{\underset{z<w,s_iz<z}{\sum}}C_z
\]
for $s_i\in{\mathcal L}(w)\{s_j |\;s_jw<w\}$. 
Since 
$C_{s_i}=q^{-\frac{1}{2}}T_i-q^{\frac{1}{2}}$, 
we can give an inductive definition of 
Kazhdan-Lusztig polynoimals as follows. 

\begin{defn}
\label{inductive definition}
\[
\begin{array}{rl}
P_{y,w}(q)&=q^{1-c}P_{s_iy,s_iw}(q)+q^cP_{y,s_iw}(q)\\
&\quad-\underset{\mu(z,s_iw)\ne0}
{\underset{y\le z\le s_iw,s_iz<z}{\sum}}
\mu(z,s_iw)q^{\frac{l(w)-l(z)}{2}}P_{y,z}(q)
\end{array}
\]
where $c=1$ if $s_iy<y$ and $c=0$ if $s_iy>y$. 
\end{defn}

That the right hand side does not depend on the 
choice of $s_i$ comes from the welldefinedness 
result. We can also construct 
$C_w$ by $C_e=1$ and 
\[
C_w=C_{ws_i}C_{s_i}-\underset{\mu(z,w)\ne0}
{\underset{z<w,zs_i<z}{\sum}}C_z
\]
for $s_i\in{\mathcal R}(w)=\{s_j |\;ws_j<w\}$, which 
leads to a similar 
inductive definition of Kazhdan-Lusztig polynomials. 
Note that ${\mathcal R}(w)={\mathcal L}(w^{-1})$. 

\begin{lemma} 
\label{properties}
(1) $P_{y,w}(0)=1$. 

\noindent
(2) $y<w, s_iy>y, s_iw<w$ imply $P_{y,w}(q)=P_{s_iy,w}(q)$. 

\noindent
(3) $P_{y^{-1},w^{-1}}(q)=P_{y,w}(q)$. 
\end{lemma}

\bigskip
(1) follows from the inductive definition 
\ref{inductive definition}. 
(2) is proved by induction on $l(y)$, which also uses 
\ref{inductive definition}. (3) follows from the first 
definition of Kazhdan-Lusztig polynomials: If we replace 
$P_{y,w}(q)$ by $P_{y^{-1},w^{-1}}(q)$ in the 
definition of $C_w$ and apply the anti-involution defined 
by $T_i\mapsto T_i$, we have $C_{w^{-1}}$. Thus 
we have the bar invariance property and the degree 
property. Hence $P_{y,w}(q)$ and $P_{y^{-1},w^{-1}}(q)$ 
must coincide. 

\bigskip
Another corollary of this construction of 
$C_w$ is the {\bf Kazhdan-Lusztig representation} of the regular 
representation, which is the matrix representation 
with respect to the basis $\{C_w\}$. 
For the left regular representation, we have
\[
T_iC_w=\left\{\begin{array}{ll}-C_w &(\mbox{if } s_iw<w)\\
\;qC_w+q^{\frac{1}{2}}C_{s_iw}+
\underset{\mu(z,w)\ne0}{\underset{z<w,s_iz<z}{\sum}}
\mu(z,w)q^{\frac{1}{2}}C_z\;&(\mbox{if }s_iw>w). 
\end{array}\right.
\]
We have the same formula for the right 
regular representation. 
We can now introduce the notion of left cells 
and right cells. 

\begin{defn}
If there exists a sequence $y=x_1,x_2,\dots,x_r=w$ 
such that ${\mathcal L}(x_i)\not\subset{\mathcal L}(x_{i+1})$, 
$\mu(x_i|x_{i+1})\ne 0$ 
for $1\le i<r$, we write $y\underset{L}{\le}w$. 

If there exists a 
sequence $y=x_1,x_2,\dots,x_r=w$ 
such that ${\mathcal R}(x_i)\not\subset{\mathcal R}(x_{i+1})$, 
$\mu(x_i|x_{i+1})\ne 0$ 
for $1\le i<r$, we write $y\underset{R}{\le}w$. 

Note that $y\underset{R}{\le}w$ if and only if 
$y^{-1}\underset{L}{\le}w^{-1}$. 
If both $y\underset{L}{\le}w$ and $w\underset{L}{\le}y$ 
hold, we write $y\underset{L}{\sim}w$. 
Similarly, if both $y\underset{R}{\le}w$ and 
$w\underset{R}{\le}y$ 
hold, we write $y\underset{R}{\sim}w$. 

These relations partition 
$S_n$ into equivalence classes, which are 
called {\bf left cells} and {\bf right cells} respectively. 
\end{defn}

At a first look, the definition of the relation 
$y\underset{L}{\le}w$ seems to be very artificial. 
To understand it in a more natural way, 
we set $q=1$ and denote $C_w|_{q=1}$ by $a(w)$. 
(The specialization 
to $q=1$ is only for simplifying the situation to 
more familiar setting of the symmetric group, and 
is not at all essential.) Then we have the following 
lemma by using Lemma \ref{properties} 
and the Kazhdan-Lusztig representation specialized 
to $q=1$. 

\begin{lemma}
\label{basal module}
Let $y\ne w$ be two elements of $S_n$. Then 
we have (1)$\Leftrightarrow$(2) where

\medskip
\noindent
(1) $s_i\in{\mathcal L}(y)\setminus{\mathcal L}(w)$ and 
$\mu(y|w)\ne 0$. 

\medskip
\noindent
(2) $a(y)$ appears in $s_ia(w)$. 
\end{lemma}

\medskip
Hence, the left regular 
representation of the symmetric group with 
the specific basis $\{a(w)\}$ gives 
a natural meaning of the relation 
$y\underset{L}{\le}w$ as follows. 

Let $\overline{V_w}^L$ be the left ideal 
uniquely defined by the following three 
conditions. 

(1)\;$a(w)\in \overline{V_w}^L$. 

(2)\;$\overline{V_w}^L$ is 
spanned by a subset of $\{a(x)\}$. 

(3)\;If a left ideal satisfies (1) and (2), 
it contains $\overline{V_w}^L$. 

\medskip
Then we have $y\underset{L}{\le}w \Leftrightarrow 
\overline{V_y}^L\subset\overline{V_w}^L$. 
Similar formula exists for $y\underset{R}{\le}w$. 

\section{RS correspondence and the left cell}
\subsection{The Kazhdan-Lusztig theorem} 
The following theorem is the theorem of Kazhdan and 
Lusztig which we are going to prove. 

\bigskip
\noindent
{\bf Theorem A}\quad 
For $y,w\in S_n$, we have 
$y\underset{L}{\sim}w \Leftrightarrow Q(y)=Q(w)$. 

\medskip
\begin{example}[The $S_3$ case] 
Left cells are $\{123\}$, $\{213,312\}$, 
$\{132,231\}$ and $\{321\}$. Their Q-symbols are 

\setlength{\unitlength}{16pt}
\begin{picture}(15,5)(-3,1)
\put(0,3){\rm1}\put(1,3){\rm2}\put(2,3){\rm3}\put(2.5,3){,}

\put(4,3.5){\rm1}\put(5,3.5){\rm3}
\put(4,2.5){\rm2}\put(5.5,3){,}

\put(7,3.5){\rm1}\put(8,3.5){\rm2}
\put(7,2.5){\rm3}\put(8.5,3){,}

\put(10,4){\rm1}
\put(10,3){\rm2}\put(10.5,3){.}
\put(10,2){\rm3}
\end{picture}

\end{example}
For the $S_4$ case, see \cite[p.20]{Shi}. 

\subsection{A theorem of Knuth}
We write $y\equiv w$ if $P(y)=P(w)$. 
To describe this equivalence relation, we 
introduce Knuth relations. 

\begin{defn}
Let $y_1\cdots y_n$ be a permutation of 
$1,\dots,n$. We set $w$ as follows. 
\begin{center}
If $y_{i+1}<y_i<y_{i+2}$, we set 
$w=y_1\cdots y_i y_{i+2} y_{i+1}\cdots y_n$. 
\end{center}
\begin{center}
If $y_{i+1}<y_{i+2}<y_i$, we set 
$w=y_1\cdots y_{i+1} y_i y_{i+2}\cdots y_n$. 
\end{center}
We have $y\equiv w$, and we 
say that $y$ and $w$ are in Knuth relation. 
\end{defn}

The following theorem is due to Knuth \cite[Theorem 6]{K}. 

\begin{thm} 
\label{Knuth's theorem}
Let $y,w\in S_n$. Then $y\equiv w$ if and only if 
these permutations are connected by a chain of 
Knuth relations. 
\end{thm}

Let $D_{ij}:=
\{\;w\;|\;ws_i<w, ws_j>w\;\}$ where $j=i\!\pm\!1$. 
If $y\in D_{ij}$, we consider the right coset 
$y<s_i,s_j>$ and take the distinguished 
coset representative $y^0$. Then we have either 
$y=y^0s_i$ or $y^0s_js_i$. We set 
$K_{ij}(y)=y^0s_is_j$ in the former case, 
and $K_{ij}(y)=y^0s_j$ in the latter case. 
Note that $K_{ij}$ is a bijective map from 
$D_{ij}$ to $D_{ji}$. 
If $j=i\!+\!1$, this is the rule to obtain $w$ 
from $y$ in the Knuth relation, and if 
$j=i\!-\!1$, this is the rule to obtain $y$ 
from $w$ in the Knuth relation. This description 
of Knuth relations is convenient for our 
purpose. The following lemma shows that 
two elements in Knuth relation are in a right cell. 

\begin{lemma} 
\label{Knuth relation}
If $w\in D_{ij}$, 
we have $K_{ij}(w)\underset{R}{\sim}w$. 
\end{lemma}
(Proof) Since $w\in D_{ij}$, we have 
$w=w^0s_i$ or $w=w^0s_js_i$ where $w^0$ is the 
distinguished coset representative of $w\langle s_i,s_j\rangle$. 
By the same proof as in Lemma \ref{basal module} 
we have $\mu(w^0s_i,w^0s_is_j)=1$, and 
$\mu(w^0s_j,w^0s_js_i)=1$. In either cases, 
we have $\mu(w|K_{ij}(w))\ne 0$. Since 
$s_i\in {\mathcal R}(w)\setminus{\mathcal R}(K_{ij}(w))$ and 
$s_j\in {\mathcal R}(K_{ij}(w))\setminus{\mathcal R}(w)$, 
we have ${\mathcal R}(w)\not\subset{\mathcal R}(K_{ij}(w))$ 
and ${\mathcal R}(w)\not\supset{\mathcal R}(K_{ij}(w))$. 
We have the result. 
\qed

\medskip
\noindent
{\bf Remark} 
We have another way to describe the 
Knuth relation as follows. 

\bigskip
Assume that $w<s_iw$ and ${\mathcal L}(w)\not\subset
{\mathcal L}(s_iw)$. 
Then we have $w^{-1}\equiv w^{-1}s_i$. 

\bigskip
In fact, if we take $s_j\in 
{\mathcal L}(w)\setminus{\mathcal L}(s_iw)$, 
the choice of $s_i,s_j$ 
leads to $s_jw<w$ and $s_js_iw>s_iw$. These 
$s_i$ and $s_j$ can not be commutable elements. 
Thus $w^{-1}\in D_{ji}$ and we are in the 
latter case in the definition of 
$K_{ij}$. If we consider the case that 
$y<s_jy$ and ${\mathcal L}(y)\not\subset
{\mathcal L}(s_jy)$, where we take $s_i\in 
{\mathcal L}(y)\setminus{\mathcal L}(s_jy)$ such that 
$y^{-1}\in D_{ij}$, 
we meet the former case in the definition 
of $K_{ij}$, and 
we have $y^{-1}\equiv y^{-1}s_j$. But 
this statement is the same as the previous 
case. 

\subsection{Preparatory results for 
the proof of Theorem A}

The following three propositions are proved in \cite{KL}. 
I avoid repetition 
as long as the readability of the proof is guaranteed. 

\begin{prop}[{\cite[Proposition 2.4]{KL}}]
\label{prop 1.2}

If $y\underset{L}{\le}w$, we have 
${\mathcal R}(y)\supset{\mathcal R}(w)$. 
\end{prop}
(Proof) It is enough to prove it for 
the case that ${\mathcal L}(y)\not\subset{\mathcal L}(w)$ 
and $\mu(y|w)\ne 0$. By Lemma \ref{basal module}, 
We have $y=s_iw>w$ for some $i$ or 
$y<w$ and $\mu(y,w)\ne 0$. 
In the former case, 
we consider the double coset $\langle s_i\rangle w
\langle s_j\rangle$ 
for each $s_j\in{\mathcal R}(w)$. 
Then we can easily conclude that $s_j\in{\mathcal R}(s_iw)$. 
Thus we have ${\mathcal R}(y)\supset{\mathcal R}(w)$. 
In the latter case, we assume to the contrary 
that there is 
$s_j\in{\mathcal R}(w)\setminus{\mathcal R}(y)$. 
By Lemma \ref{properties}(3), our assumption 
$\mu(y|w)\ne 0$ is equal to 
$\mu(y^{-1}|w^{-1})\ne 0$. We also have 
$s_j\in{\mathcal L}(w^{-1})\setminus{\mathcal L}(y^{-1})$. 
By Lemma \ref{basal module}, we know 
that $a(w^{-1})$ appears in $s_ja(y^{-1})$. 
Since $y<w$, we have $w^{-1}>y^{-1}$ and 
thus we have $w^{-1}=s_jy^{-1}$. 
By the same argument in the former case, 
$w=ys_j>y$ implies ${\mathcal L}(y)\subset{\mathcal L}(w)$. 
It contradicts our assumption that 
${\mathcal L}(y)\not\subset{\mathcal L}(w)$. 
\qed

\begin{prop}[{\cite[Theorem 4.2]{KL}}] 
\label{prop 1.3}

If $y\ne w\in D_{ij}$ and $\mu(y|w)\ne 0$, then 
$\mu(K_{ij}(y)|K_{ij}(w))\ne 0$. 
\end{prop}

\bigskip
\noindent
{\bf Remark} 
By the definition of $D_{ij}$, there are 
two possibilities for $y$ and $w$ respectively. 
Namely, 
\[
\begin{array}{rl}
ys_i<y\!\!&=K_{ij}(y)s_j<ys_j=K_{ij}(y)<K_{ij}(y)s_i,\\
ys_j>y\!\!&=K_{ij}(y)s_i>ys_i=K_{ij}(y)>K_{ij}(y)s_j,\\
&{\mbox or}\\
ws_i<w\!\!&=K_{ij}(w)s_j<ws_j=K_{ij}(w)<K_{ij}(w)s_i,\\
ws_j>w\!\!&=K_{ij}(w)s_i>ws_i=K_{ij}(w)>K_{ij}(w)s_j.
\end{array}
\]
Let $y_i,w_i,\;(i=1,2)$ and $s,t$ be as follows. 

\noindent
(a) If both $y$ and $w$ are in the former case, 
we set 
\[
y_1=K_{ij}(y),\;y_2=y,\;s=s_j,\;t=s_i,\;
w_1=K_{ij}(w),\;w_2=w.
\]

\noindent
(b) If $w$ is in the latter case, 
we set 
\[
y_1=y,\;y_2=K_{ij}(y),\;s=s_i,\;t=s_j,\;
w_1=w,\;w_2=K_{ij}(w).
\]

\noindent
(c) If $y$ is in the latter case and 
$w$ is in the former case, we set 
\[
w_1=K_{ij}(y),\;w_2=y,\;s=s_j,\;t=s_i,\;
y_1=K_{ij}(w),\;y_2=w.
\]
Then we are reduced to the following two cases. 
\[
\begin{array}{rl}
(1)\;y_2t<y_2=y_1s<y_1<y_1t,&\;w_2t<w_2=w_1s<w_1<w_1t,\\
(2)\;y_2t<y_2=y_1s<y_1<y_1t,&\;w_1s<w_1<w_1t=w_2<w_2s.
\end{array}
\]
We have to show $\mu(y_1|w_1)=\mu(y_2|w_2)$ for these 
two cases. Then we have come to the beginning 
of the proof in \cite[Theorem 4.2(i\!i\!i)]{KL}. 

\bigskip

\begin{prop}[{\cite[Corollary 4.3]{KL}}] 
\label{Knuth move}

Let $y, w\in D_{ij}$. Then $y\underset{L}{\sim}w$ implies 
$K_{ij}(y)\underset{L}{\sim}K_{ij}(w)$. 
\end{prop}
(Proof) We can assume that ${\mathcal L}(y)\not\subset{\mathcal L}(w)$, 
${\mathcal L}(y)\not\supset{\mathcal L}(w)$, and $\mu(y|w)\ne 0$. 
By the same proof as in Proposition \ref{prop 1.2}, 
we have ${\mathcal L}(ys_i)={\mathcal L}(y)$. Hence we have 
${\mathcal L}(K_{ij}(y))\not\subset{\mathcal L}(K_{ij}(w))$ and 
${\mathcal L}(K_{ij}(y))\not\supset{\mathcal L}(K_{ij}(w))$. 
We also have $\mu(K_{ij}(y)|K_{ij}(w))\ne 0$ 
by Proposition \ref{prop 1.3}. We are through. 
\qed

\subsection{Proof of Theorem A}

One implication is easy. 

\begin{prop}
\label{one implication}
If $Q(y)=Q(w)$, then $y\underset{L}{\sim}w$. 
\end{prop}
(Proof) Since $Q(y)=P(y^{-1})$ and $Q(w)=P(w^{-1})$, $y^{-1}$ 
is connected to $w^{-1}$ by a chain of Knuth relations. 
Thus it is enough to prove that $w^{-1}=K_{ij}(y^{-1})$ 
$(y^{-1}\in D_{ij})$ 
implies $y\underset{L}{\sim}w$. But Lemma \ref{Knuth relation} 
shows that $y^{-1}\underset{R}{\sim}w^{-1}$, 
which is $y\underset{L}{\sim}w$. 
\qed

\bigskip
It remains to prove that $y\underset{L}{\sim}w$ implies 
$Q(y)=Q(w)$. For each partition $\pi$, we define 
a standard tableau $P_\pi$ by setting the entries of 
the $i$ th column of $P_\pi$ 
to be $\sum_{j=1}^{i-1}l_j+1,\dots,\sum_{j=1}^i l_j$ from top to 
bottom, where $l_1,l_2,\dots$ are column lengths of $\pi$. 
We denote the shapes of $Q(y)$ and $Q(w)$ by $\pi_1$ and $\pi_2$ 
respectively. We define ${\hat y},$ ${\hat w}$ by 
$(P({\hat y}),Q({\hat y}))=(P_{\pi_1},Q(y))$ and 
$(P({\hat w}),Q({\hat w}))=(P_{\pi_2},Q(w))$. 
By Proposition \ref{one implication}, we have 
$y\underset{L}{\sim}{\hat y}$ and $w\underset{L}{\sim}{\hat w}$. 
Thus we have ${\hat y}\underset{L}{\sim}{\hat w}$. To prove that 
$Q({\hat y})=Q({\hat w})$, we define $y'$ and $w"$ by 
$(P(y'),Q(y'))=(P_{\pi_1},P_{\pi_1})$ and 
$(P(w"),Q(w"))=(P_{\pi_2},P_{\pi_2})$. By the theorem of Knuth, 
we can write 
\[
\begin{array}{rl}
y'&=K_{i_1j_1}\circ\cdots\circ K_{i_rj_r}({\hat y})\\
w"&=K_{i_1'j_1'}\circ\cdots\circ K_{i_s'j_s'}({\hat w})
\end{array}
\]
We shall define $w'$ and $y"$ by 
\[
\begin{array}{rl}
w'&=K_{i_1j_1}\circ\cdots\circ K_{i_rj_r}({\hat w})\\
y"&=K_{i_1'j_1'}\circ\cdots\circ K_{i_s'j_s'}({\hat y})
\end{array}
\]
Recall that Proposition \ref{prop 1.2} tells that 
${\mathcal R}({\hat y})={\mathcal R}({\hat w})$. Hence 
${\hat y}\in D_{i_rj_r}$ implies ${\hat w}\in D_{i_rj_r}$. 
We then have welldefined $K_{i_rj_r}({\hat w})$, which satisfies 
$K_{i_rj_r}({\hat y})\underset{L}{\sim}
K_{i_rj_r}({\hat w})$ by Proposition 
\ref{Knuth move}. We continue the argument and conclude that 
these $y"$ and $w'$ are welldefined. 
$y'$ and $w'$ satisfy ${\mathcal R}(y')={\mathcal R}(w')$, 
$P(y')=Q(y')=P_{\pi_1}$ and $P(w')=P_{\pi_2}$. 
Similarly, $y"$ and $w"$ satisfy 
${\mathcal R}(y")={\mathcal R}(w")$, $P(y")=P_{\pi_1}$ and 
$P(w")=Q(w")=P_{\pi_2}$. 
Note that $y'$ is the permutation 
\[
l_1,l_1\!-\!1,\cdots,1,l_1\!+\!l_2,\cdots,l_1\!+\!1,\cdots. 
\]
Similarly, $w"$ is the permutation 
\[
l_1',l_1'\!-\!1,\cdots,1,l_1'\!+\!l_2',\cdots,l_1'\!+\!1,\cdots. 
\]
where we denote column lengths of $\pi_2$ by 
$l_1',l_2',\dots$. 
Since ${\mathcal R}(y')={\mathcal R}(w')$, the first $l_1$ letters of 
$w'$ are in the decreasing order, the next $l_2$ letters are 
in the decreasing order, etc. Similarly, the first $l_1'$ 
letters of $y"$ are in the decreasing order, the next $l_2'$ 
letters are in the decreasing order, etc. 

By inserting the first $l_1$ letters of $w'$ to $\emptyset$, 
we know that the first column of 
$\pi_2$ must have the length equal or greater than $l_1$. 
By using $y"$, we have the opposite 
inequality. We have $l_1=l_1'$. It also implies that 
the next $l_2$ decreasing letters of $w'$ do not produce bumping, 
since if otherwise we have $l_1'>l_1$. 
Thus we have that $l_2'\ge l_2$. We use $y"$ to have the opposite 
inequality. Continuing the same argument, we conclude that 
$\pi_1=\pi_2$ and $Q(w')=P_{\pi_2}$. (We also have 
$Q(y")=P_{\pi_1}$.) Therefore, we have $y'=w'$, which 
implies $\hat{y}=\hat{w}$. We have proved $Q(y)=Q(w)$. 
\qed

\subsection{Theorem A in the crystal base theory context} 

An occurence of the Robinson-Schensted algorithm in the 
tensor product representation of the vector representation 
of $U_q({\mathfrak gl}_n)$ was first observed in \cite{DJM}. 
The tensor product representation itself can be viewed as 
an example of Demazure modules \cite[Theorem 3.1]{KMOTU}, 
and we may consider generalization into this direction, 
but we restrict ourselves to the original case. Then 
the crystal base is induced by the canonical base and 
we now have a good understanding of the base (see \cite{SV}) and 
of the Robinson-Schensted algorithm in the crystal base 
theory context. 

Let $U_q$ be the quantum algebra of ${\mathfrak gl}_r$, and 
$\Delta$ be Lusztig's coproduct:
\[
\begin{array}{rl}
\Delta(e_i)=& e_i\otimes 1 + 
q^{\epsilon_i-\epsilon_{i\!+\!1}}\otimes e_i
\quad (i=1,\dots,r\!\!-\!\!1),\\
\Delta(f_i)=& 1 \otimes f_i + f_i\otimes 
q^{-\epsilon_i+\epsilon_{i\!+\!1}}\quad (i=1,\dots,r\!\!-\!\!1),\\
\Delta(q^h)=& q^h\otimes q^h \quad (\;h\in{\mathbb Z}\epsilon_1+
\cdots+{\mathbb Z}\epsilon_n).
\end{array}
\]
Let $V={\mathbb Q}(q)^r$ be its vector representation given by 
\[
e_i=E_{i,i+1},\quad f_i=E_{i+1,i},\quad q^{\epsilon_i}
=qE_{ii}+\sum_{j\ne i}E_{jj}
\]
where $E_{ij}$ are matrix units. 
Natural base elements $v_1=(1,0,\dots,0)^{\rm T}, 
v_2=(0,1,\dots,0)^{\rm T},\dots$ induce a base at $q=\infty$ 
in the sense of Kashiwara-Lusztig. In the following, 
we exclusively work 
with bases at infinity, and call them 
crystal bases instead of bases at $q=\infty$. 
We set 
$L=\oplus_{i=1}^r {\mathbb Q}[q^{-1}]_{(q^{-1})}v_i$, 
$B=\{\;v_i\;{\rm mod}\;q^{-1}L\}\subset L/q^{-1}L$. 
Then $(L,B)$ is the crystal base of $V$ stated above, and 
$V^{\otimes n}$ has $(\;L^{\otimes n},B^{\otimes n})$ 
as its crystal base. To describe the tensor structure on 
$B^{\otimes n}$, we introduce $\varphi_i(b), \epsilon_i(b)$ by 
\[
\varphi_i(b)=\mbox{\rm max}\{\;k |\; \tilde f_i^k(b)\ne 0\}, \quad
\epsilon_i(b)=\mbox{\rm max}\{\;k |\; \tilde e_i^k(b)\ne 0\}
\]
where $\tilde e_i$ and $\tilde f_i$ are Kashiwara operators. 
Then 
\[
\begin{array}{rl}
\tilde e_i(b_1\otimes b_2)&=\left\{\begin{array}{rl}
b_1\otimes\tilde e_i(b_2)\quad&(\epsilon_i(b_1)\le\varphi_i(b_2)),\\
\tilde e_i(b_1)\otimes b_2\quad&(\epsilon_i(b_1)>\varphi_i(b_2)),
\end{array}\right.
\\
\tilde f_i(b_1\otimes b_2)&=\left\{\begin{array}{rl}
b_1\otimes \tilde f_i(b_2)\quad&(\epsilon_i(b_1)<\varphi_i(b_2)),\\
\tilde f_i(b_1)\otimes b_2\quad&(\epsilon_i(b_1)\ge\varphi_i(b_2)).
\end{array}\right.
\end{array}
\]
Let $V_q(\lambda)$ be the irreducible highest module of $U_q$ 
associated with $\lambda=\sum \lambda_i\epsilon_i$. We identify 
$\lambda$ with the corresponding Young diagram. Then it is well known 
that $V_q(\lambda)\otimes V$ is multiplicity free. Hence, we can 
uniquely define the submodule of $V^{\otimes n}$ 
for each increasing sequence of Young diagrams. 
We identify the increasing sequence with a standard tableau 
$Q$, which we call the recording tableau. We denote the 
submodule by $V_q(Q)$. If the shape of $Q$ is $\lambda$, 
we have $V_q(Q)\simeq V_q(\lambda)$. 

\begin{prop}

(1) Let $Q$ be a standard tableau of size $n\!\!-\!\!1$, and 
${\mathcal T}_Q$ be the set of tableaux obtained from $Q$ by adding 
$\boxed{n}$. Let $(L(Q),B(Q))$ be a crystal base of $V_q(Q)$. 
We set $L(T)=(L(Q)\otimes L)\cap V_q(T)$. Then we have 
\[
L(Q)\otimes L=\bigoplus_{T\in{\mathcal T}_Q} L(T). 
\]
We nextly set $B(T)=(B(Q)\otimes B)\cap(L(T)/q^{-1}L(T))$. 
Then we have 
\[
(L(Q)\otimes L,B(Q)\otimes B)=\bigoplus_{T\in{\mathcal T}_Q} (L(T),B(T)). 
\]

\noindent
(2) Let $L(Q)=L^{\otimes n}\cap V_q(Q)$. Then we have 
$L^{\otimes n}=\oplus L(Q)$. If we further set 
$B(Q)=B^{\otimes n}\cap(L(Q)/q^{-1}L(Q))$, we have 
$(L^{\otimes n},B^{\otimes n})=\oplus (L(Q),B(Q))$. 
\end{prop}
(Proof) (1) Let $v_T$ be the highest weight vector which 
generates the highest weight space of $L(T)$. 
Since $L(T)$ and the lattice generated by 
$\tilde f_{i_1}\cdots\tilde f_{i_*}v_T$ are crystal lattices of 
$V_q(T)$, the uniqueness theorem of crystal bases concludes that 
they coincide. The uniqueness theorem also guarantees that 
there exists an automorphism of $V_q(Q)\otimes V$ such that 
it maps $L(Q)\otimes L$ to $\oplus L(T)$. 
Since $V_q(Q)\otimes V$ is multiplicity free, the automorphism 
is scalar multiplication on each $V_q(T)$. Thus by looking at 
highest weight spaces, we have that the automorphism is the 
identity. By descending induction on weights, we can prove 
$B(Q)\otimes B=\sqcup B(T)$. 

\noindent
(2) We prove it by induction on $n$. Assume that it holds for 
$n$. Then we have 
$(L(Q)\otimes L,B(Q)\otimes B)=\oplus (L(T),B(T))$ by (1) where 
\[
\begin{array}{rl}
L(T)&=(L(Q)\otimes L)\cap V_q(T)
=\left((L^{\otimes n}\cap V_q(Q))\otimes L\right)\cap V_q(T)\\
&=\left(L^{\otimes n\!+\!1}\cap V_q(Q)\otimes L\right)\cap V_q(T)\\
&=\left(L^{\otimes n\!+\!1}\cap V_q(Q)\otimes V\right)\cap V_q(T)\\
&=L^{\otimes n\!+\!1}\cap V_q(T).
\end{array}
\]
\qed

\bigskip
The crystal graph of $V_q(\lambda)$ has description 
in terms of semistandard tableaux as follows \cite{KN}. 

We write $\boxed{i}$ for $v_i\;\mbox{\rm mod\;}q^{-1}\in B$. 
Let $B(\lambda)$ be the set of column strict semi-standard 
tableaux of shape $\lambda$. For each $T\in B(\lambda)$, 
we read its entries row by row, starting from the bottom row. 
This reading gives an injection from $B(\lambda)$ to 
$B^{\otimes n}$. 
For example, we have 
\[
\begin{array}{cccc}
1 & 1 & 2 & 4 \\
2 & 3 &   &   \\
4 &   &   &   \end{array}\qquad \Longrightarrow 
\boxed{4}\otimes\boxed{2}\otimes\boxed{3}\otimes\boxed{1}
\otimes\boxed{1}\otimes\boxed{2}\otimes\boxed{4}. 
\]
We induce the crystal structure on $B(\lambda)$ 
through this inclusion: the Kashiwara 
operators $\tilde e_i, \tilde f_i$ 
act on these monomial tensors by changing the leftmost 
$\boxed{i\!\!+\!\!1}$ 
or the rightmost $\boxed{i}$ of the sequence which is 
obtained by 
removing consecutive $\boxed{i\!\!+\!\!1}\otimes\boxed{i}$ as many 
as possible. Thus $B(\lambda)$ is stable by 
Kashiwara operators. This embedding is in fact 
the embedding of $B(\lambda)$ into the set of permutation 
tableaux by {\it jeu de taquin} moves, and the inverse is 
given by taking $P$-symbols, namely by the insertion algorithm. 
See \cite{Fo3} for example. 

Let $Q$ be a standard tableau of shape $\lambda$. We 
identify $B(Q)$ with $B(\lambda)$. Note that there exists 
a unique isomorphism of the crystals $(L(Q),B(Q))$ and 
$(L(\lambda),B(\lambda))$. The following is the modern version 
of the Date-Jimbo-Miwa theorem. 
We refer \cite{BKK} for its generalization to super algebras. 

\begin{thm} 
We identify $B(Q)$ with $B(\lambda)$ as above. Then the following 
hold. 

\bigskip
\noindent
(1) If $b=\boxed{i_1}\otimes \cdots\otimes\boxed{i_n}\in B(Q)$, 
then the Q-symbol of 
$\emptyset\leftarrow i_1i_2\cdots i_n$ 
is $Q$. 

\noindent
(2) Let $P(b)$ be the P-symbol of 
$\emptyset\leftarrow i_1i_2\cdots i_n$. 
Then the identification of $B(Q)$ with $B(\lambda)$ is given by 
the map $b \mapsto P(b)$. 
\end{thm}
(Proof) We first recall that 
the bumping procedure $(T,i) \mapsto (T\!\leftarrow\! i)$ gives 
the isomorphism of crystals between $B(\lambda)\otimes B$ and 
$\sqcup_{|\mu/\lambda|=1} B(\mu)$. (As I have explained, 
we can think of the insertion via {\it jeu de taquin} moves. 
Hence it is enough to establish the 
isomorphism for a {\it jeu de taquin} move, which is easy.) 

This isomorphism leads to a crystal automorphism 
on $B(Q)\otimes B$ as folows.
\[
B(Q)\otimes B\stackrel{\sim}{\rightarrow} B(\lambda)\otimes B 
\stackrel{\sim}{\rightarrow}\bigsqcup_{|\mu/\lambda|=1}B(\mu)
\stackrel{\sim}{\leftarrow}
\bigsqcup_{T\in{\mathcal T}_Q} B(T)=B(Q)\otimes B
\]
where the second isomorphism is given by the insertion algorithm. 
Since $V_q(Q)\otimes V$ is multiplicity free, the automorphism 
must be the identity. Hence the isomorphim 
$B(T)\stackrel{\sim}{\rightarrow}B(\mu)$ for $T\in{\mathcal T}_Q$ of shape $\mu$ 
is given by restricting the following isomorphism to $B(T)$. 
Note again that the second isomorphism is given by the insertion 
algorithm. 
\[
B(Q)\otimes B \stackrel{\sim}{\rightarrow}
B(\lambda)\otimes B \stackrel{\sim}{\rightarrow} 
\bigsqcup_{|\mu/\lambda|=1} B(\mu)
\]
Thus if the Robinson-Schensted algorithm gives 
the isomorphism 
$B(Q)\stackrel{\sim}{\rightarrow}B(\lambda)$ 
such that the Q-symbols of the elements in its image are 
constant $Q$, then the Robinson-Schensted algorithm 
gives the isomorphim $B(T)\stackrel{\sim}{\rightarrow}B(\mu)$, 
and the Q-symbols of the elements in its image are 
constant $T$. Therefore the induction proceeds. \qed

\bigskip
We now turn to the $q^2$-Schur algebra. We refer \cite{Du} 
for the details. 
We consider the Hecke algebra whose deformation parameter 
$q$ is replaced by $q^2$. We also denote it by 
${\mathcal H}_n$ by abuse of notion. 
$V^{\otimes n}$ has ${\mathcal H}_n$ action given by 
\[
v_{i_1}\otimes\cdots\otimes v_{i_n}T_k=
\left\{\begin{array}{ll}
qv_{i_1}\otimes\cdots\otimes v_{i_{k+1}}\otimes 
v_{i_k}\otimes\cdots\otimes v_{i_n}\quad&(i_k>i_{k+1})\\
q^2v_{i_1}\otimes\cdots\otimes v_{i_n}\quad&(i_k=i_{k+1})\\
qv_{i_1}\otimes\cdots\otimes v_{i_{k+1}}\otimes 
v_{i_k}\otimes\cdots\otimes v_{i_n}&\\
\qquad +(q^2-1)v_{i_1}\otimes\cdots\otimes v_{i_n}
\quad&(i_k<i_{k+1})\end{array}\right.
\]
It commutes with $U_q$ action. The 
endomorphism ring ${\rm End}_{{\mathcal H}_n}(V^{\otimes n})$ 
is called the {\bf $q$-Schur algebra}, which is denoted 
by ${\mathcal S}_{r,n}$. It is well known that 
it is a quotient algebra of $U_q$. 
If we denote the $\mu$-weight space of $V^{\otimes n}$ by 
$V_\mu$, then we obviously have 
${\mathcal S}_{r,n}=\oplus_{\mu,\nu}{\rm End}_{{\mathcal H}_n}
(V_\nu,V_\mu)$. 

We now assume $r=n$ and set 
$\omega=\epsilon_1\!+\!\cdots\!+\!\epsilon_n$. 
Then ${\mathcal H}_n\simeq {\rm End}_{{\mathcal H}_n}
(V_\omega,V_\omega)$, and we can identify ${\mathcal H}_n$ 
with the subalgebra of ${\mathcal S}_{n,n}$. 

On the other hand, if we set $x_\mu=\sum_{w\in S_\mu} T_w$ 
where $S_\mu$ is the Young subgroup associated with $\mu$, 
the weight space $V_\mu$ is isomorphic to $x_\mu{\mathcal H}_n$. 
Hence we can also identify $V_\omega$ with ${\mathcal H}_n$. 
This identification is given by 
\[
v_{w_1}\otimes\cdots\otimes v_{w_n}\mapsto 
(q^2)^{-l(ww_0)/2}T_{ww_0}. 
\]
In particular, the Kazhdan-Lusztig basis element $C_w'$ is identified with 
\[
\sum_{y\le w} P_{y,w}(q^2)q^{l(y)-l(w)}v_{y_n}\otimes\cdots
\otimes v_{y_1}. 
\]
We have 
$P_{y,w}(q^2)q^{l(y)-l(w)}\in{\mathbb Z}[q^{-1}]$, 
and $C_w'\equiv v_{w_n}\otimes\cdots\otimes v_{w_1}\;
{\rm mod\;}q^{-1}$. 

\bigskip
The tensor space and the $q^2$-Schur algebra have bar operations, 
which satisfy $\bar{x}\bar{v}=\overline{xv}\;
(x\in{\mathcal S}_{n,n}, v\in V^{\otimes n})$, and 
$\overline{v_n\otimes\cdots\otimes v_1}=v_n\otimes\cdots\otimes v_1$. 
The bar operation on the tensor space coincides 
with the bar operation introduced in 2.2 if restricted 
to ${\mathcal H}=V_\omega\subset V^{\otimes n}$. 

By these reasons, we conclude that
these are canonical basis elements arising from 
the crystal base we have considered above. We also remark that 
the canonical basis of the $q^2$-Schur algebra is 
the image of the canonical basis 
of the modified quantized enveloping algebra by the work 
\cite{SV}. 
In fact, because of $(V^{\otimes n})_\omega=\oplus S_q(Q)$ where 
$S_q(Q)=V_q(Q)_\omega$, these Kazhdan-Lusztig basis 
elements are partitioned into the disjoint union 
$\sqcup B(Q)_\omega$ at $q=\infty$. 

\bigskip
Recall that these $C_w'$ are obtained from $C_w$ by applying 
a ${\mathbb Q}$-algebra automorphism of ${\mathcal H}_n$. Thus 
the vector spaces 
$S_{\underset{L}{\le}w}$, $S_{\underset{L}{<}w}$ 
generated by $\{C_y'|y\underset{L}{\le}w\}$, 
$\{C_y'|y\underset{L}{<}w\}$ respectively 
are ${\mathcal H}_n$-modules. It is known that the factor module 
$S_{\underset{L}{\le}w}/S_{\underset{L}{<}w}$ is irreducible. 
We now take the $U_q$-submodules 
$V_{\underset{L}{\le}w}$, $V_{\underset{L}{<}w}$ 
of $V^{\otimes n}$ generated by 
$S_{\underset{L}{\le}w}$, $S_{\underset{L}{<}w}$ respectively. 
By applying compositions of $\tilde e_i, \tilde f_i$ to 
$\{C_y'|y\underset{L}{\le}w\}$, $\{C_y'|y\underset{L}{<}w\}$, 
we also have crystal bases 
of $V_{\underset{L}{\le}w}$ and $V_{\underset{L}{<}w}$, 
which we denote by $(L_{\underset{L}{\le}w},B_{\underset{L}{\le}w})$, 
$(L_{\underset{L}{<}w},B_{\underset{L}{<}w})$. 
$B_{\underset{L}{<}w}$ is a union of connected components of 
$B_{\underset{L}{\le}w}$. 
Since $V_{\underset{L}{\le}w}/V_{\underset{L}{<}w}$ is irreducible, 
$B_{\underset{L}{\le}w}\setminus B_{\underset{L}{<}w}$ 
coincides with one of $B(Q)_\omega$. 
Therefore, we have Theorem A again. 

\subsection{Theorem A derived from the primitive ideal theory} 

\begin{defn} 
The annihilator ideal of $L(\lambda)$ in $U({\mathfrak g})$ 
is denoted by $I(\lambda):=Ann(L(\lambda))$, and 
is called a {\bf primitive ideal}. 
\end{defn}

The following is a theorem of Joseph.

\bigskip
\noindent
{\bf Theorem B}\quad 
$Q(y)=Q(w) \Leftrightarrow I(y\cdot0)=I(w\cdot0)$. 

\bigskip
By the translation principle, 
$0$ can be replaced by any dominant integral. 

The proof of this theorem depends on the following proposition. 

\begin{prop}
\label{fundamental lemma}
\begin{description}
\item[(1)] 
Let $y, w\in D_{ij}$ and assume that $I(y\cdot0)\subset 
I(w\cdot0)$, then we have $I(K_{ij}(y)\cdot0)\subset
I(K_{ij}(w)\cdot0)$. 
\item[(2)]
If $Q(y)=Q(w)$, we have $I(y\cdot0)=I(w\cdot0)$. 
\end{description}
\end{prop}

\bigskip
(1) is proved in \cite[Satz 5.9]{Ja1}. (2) is proved in 
\cite[Satz 5.18]{Ja1}. Once this proposition is established, 
the proof of Theorem B goes precisely the same as the proof 
of Theorem A. 

\bigskip
By \cite[Corollar 6.26]{Ja1}, \cite[Satz 7.9]{Ja1}, 
\cite[Satz 7.12]{Ja1}, we have the following theorem of Vogan. 
We state it in weaker form since it is enough for our purpose. 

\bigskip
\noindent
{\bf Theorem C}\quad 
Let $\lambda,\mu_1,\mu_2$ be dominant integral weights. 
Then we have that $I(y\cdot\lambda)\subset I(w\cdot\lambda)$ 
holds if and only if there exists a finite dimensional module $E$ 
such that 
\[
\left[L(y^{-1}\cdot\mu_1)\otimes E:L(w^{-1}\cdot\mu_2)\right]\ne0. 
\]

\bigskip
This theorem leads to Theorem D below. Recall that 
Kazhdan-Lusztig conjecture states that if we define $a(y,w)$ by 
\[
L(y\cdot0)=\sum_{y\le w} a(y,w)M(w\cdot0), 
\]
we have $a(y,w)=(-1)^{l(w)-l(y)}P_{w_0w,w_0y}(1)$. 
This is proved by Brylinski and Kashiwara, Beilinson and Bernstein. 
Thus, there is a linear isomorphism between 
$K_0({\mathcal O}_0)$ and ${\mathbb Z}W$ which sends 
$M(w_0w^{-1}\cdot0)$ to $w$ and $L(w_0w^{-1}\cdot0)$ to $a(w)$. 
By introducing $W$-action on $K_0({\mathcal O}_0)$ by 
$\tau M(w_0w^{-1}\cdot0)=M(w_0w^{-1}\tau^{-1}\cdot0)$, we can make 
it into a $W$-module isomorphism. Hence it is possible to 
translate statements for $K_0({\mathcal O}_0)$ to those for 
the Weyl group. 
The following theorem is due to Joseph and Vogan. 
The formulation is due to Joseph \cite{Jo5}, and 
Vogan gives the proof in proving Theorem C. 
See \cite[Lemma 14.9]{Ja1} for the proof. 

\bigskip
\noindent
{\bf Theorem D}\quad 
$I(yw_0\cdot0)\subset I(ww_0\cdot0) \Leftrightarrow 
a(w)\in \overline{V_y}^L$. 

\bigskip
This theorem shows that $y\underset{L}{\sim}w \Leftrightarrow 
I(yw_0\cdot0)=I(ww_0\cdot0)$. We then use Theorem B 
to conclude that $y\underset{L}{\sim}w \Leftrightarrow 
Q(yw_0)=Q(ww_0)$. Sch\"utzenberger's theorem \cite{Sch1} tells that 
if we apply evacuation procedure to $Q(w)$, we obtain 
the transpose of $Q(ww_0)$ \cite[Theorem 3.114]{S}. Thus 
we have established Theorem A again.

\bigskip
\noindent
Tokyo University of Mercantile Marine,

\noindent
Etchujima 2-1-6, Koto-ku, Tokyo 135-8533, Japan 

\noindent
ariki@ipc.tosho-u.ac.jp

\end{document}